\documentclass{article}

\title{\LARGE \textbf{A Common Generalization of Dirac's two Theorems}}
\author{Zh.G. Nikoghosyan}

\begin{document}

\maketitle

\begin{abstract}
Let $G$ be a 2-connected graph of order $n$ and let $c$ be the circumference - the order of a longest cycle in $G$.  In this paper we present a sharp lower bound for the circumference based on minimum degree $\delta$ and $p$ - the order of a longest path in $G$. This is a common generalization of two earlier classical results for 2-connected graphs due to Dirac: (i) $c\ge \min\{n,2\delta\}$; and (ii) $c\ge\sqrt{2p}$. Moreover, the result is stronger than (ii). \\


\end{abstract}

\section{Introduction}

We consider only undirected graphs with no loops or multiple edges. Let $G$ be a graph of order $n$ and let $c$ be the circumference - the order of a longest cycle in $G$. A graph $G$ is hamiltonian if $G$ contains a Hamilton cycle, that is a simple cycle $C$ with $|C|=c=n$. A good reference for any undefined terms is \cite{[1]}.

The earliest two nontrivial lower bounds for the circumference were developed in 1952 due to Dirac \cite{[3]} in terms of minimum degree $\delta$ and $p$ - the order of a longest path in $G$, respectively. \\

\noindent\textbf{Theorem A} \cite{[3]}. In every 2-connected graph, $c\ge \min\{n,2\delta\}$. \\

\noindent\textbf{Theorem B} \cite{[3]}. In every 2-connected graph, $c\ge\sqrt{2p}$.\\

In this paper we present a common generalization of Theorem A and Theorem B, including both $\delta$ and $p$ in a common relation as parameters.\\

\noindent\textbf{Theorem 1}. Let $G$ be a 2-connected graph. Then
$$
c\ge\left\{
\begin{array}{lll}
p & \mbox{when} & \mbox{ } p\le 2\delta,\\  \\
 p-1 & \mbox{when} & \mbox{ }2\delta+1\le p\le 3\delta-2,\\  \\ \sqrt{2p-10+(\delta-\frac{7}{2})^2}+\delta+\frac{1}{2} & \mbox{when} & \mbox{ }%
p\ge3\delta-1.
\end{array}
\right.
$$

 Since $G$ is 2-connected, we can assume that $n\ge 3$. If  $p\le 2\delta$ then by Theorem 1, $c\ge p$, implying that $c=p=n$ ($G$ is hamiltonian) and $c=p> \sqrt{2p}$. Next, if $2\delta+1\le p\le 3\delta-2$ then by Theorem 1, $c\ge p-1$. Since $p\ge 5$, we have $c\ge p-1\ge 2\delta$ and $c\ge p-1> \sqrt{2p}$. Finally, if $p\ge 3\delta-1$ then
 $$
 \sqrt{2p-10+\left(\delta-\frac{7}{2}\right)^2}\ge \sqrt{2(3\delta-1)-10+\left(\delta-\frac{7}{2}\right)^2}=\delta-\frac{1}{2},
 $$
 implying that
 $$
 \left(\delta+\frac{1}{2}\right)\sqrt{2p-10+\left(\delta-\frac{7}{2}\right)^2}+\delta^2-3\delta+\frac{5}{4}
 $$
 $$
 \ge\left(\delta+\frac{1}{2}\right)\left(\delta-\frac{1}{2}\right)+\delta^2-3\delta+\frac{5}{4}=(\delta-1)(2\delta-1)>0.
 $$
 Observing that the inequality
 $$
 \left(\delta+\frac{1}{2}\right)\sqrt{2p-10+\left(\delta-\frac{7}{2}\right)^2}+\delta^2-3\delta+\frac{5}{4}>0
 $$
 is equivalent to
 $$
 \sqrt{2p-10+\left(\delta-\frac{7}{2}\right)^2}+\delta+\frac{1}{2}>\sqrt{2p},
 $$
we conclude (by Theorem 1) that $c>\sqrt{2p}$.

Thus, Theorem 1 is not weaker than Theorem A and is stronger than Theorem B.

To show that Theorem 1 is best possible in a sense, observe first that in general, $p\ge c$, that is $c=p$ when $p\le 2\delta$, implying that the bound $c\ge p$ in Theorem 1  cannot be replaced by $c\ge p+1$. On the other hand, the graph $K_{\delta,\delta+1}$ with $p=2\delta+1$ and $c=2\delta=p-1$ shows that the condition $p\le 2\delta$ cannot be relaxed to $p\le 2\delta+1$. In addition, the graph $K_{\delta,\delta+1}$ with $c=p$ shows that the bound $c\ge p-1$ (when $2\delta+1\le p \le 3\delta-2$) cannot be replaced by $c\ge p$. Further, the graph $K_2+3K_{\delta-1}$ with $n=p=3\delta-1$ and $c=2\delta\le p-2$ shows that the condition $p\le 3\delta-2$ cannot be relaxed to $p\le 3\delta-1$. Finally, the same graph $K_2+3K_{\delta-1}$ with $p=3\delta-1$ and
$$
c= 2\delta= \sqrt{2p-10+\left(\delta-\frac{7}{2}\right)}+\delta+\frac{1}{2}
$$
shows that the bound $\sqrt{2p-10+(\delta-\frac{7}{2})}+\delta+\frac{1}{2}$ in Theorem 1 cannot be improved to $\sqrt{2p-10+(\delta-\frac{7}{2})}+\delta+1$.

To prove Theorem 1, we use the result of Ozeki and Yamashita \cite{[4]} for the case when $2\delta+1\le p\le 3\delta-2$.\\

\noindent\textbf{Theorem C} \cite{[4]}. Let $G$ be a 2-connected graph. Then either $(i)$ $c\ge p-1$ or $(ii)$ $c\ge 3\delta-3$ or $(iii)$ $\kappa=2$ and $p\ge 3\delta-1$.

\section{Notation and preliminaries}

The set of vertices of a graph $G$ is denoted by $V(G)$ and the set of edges by $E(G)$.  The neighborhood of a vertex $x\in V(G)$ will be denoted by $N(x)$. We use  $d(x)$ to denote $|N(x)|$.

Paths and cycles in a graph $G$ are considered as subgraphs of $G$. If $Q$ is a path or a cycle, then the order of $Q$, denoted by $|Q|$, is $|V(Q)|$. We write a cycle $Q$ with a given orientation by $\overrightarrow{Q}$. For $x,y\in V(Q)$, we denote by $x\overrightarrow{Q}y$ the subpath of $Q$ in the chosen direction from $x$ to $y$. For $x\in V(Q)$, we denote the successor and the predecessor of $x$ on $\overrightarrow{Q}$ by $x^+$ and $x^-$, respectively. We say that vertex $z_1$ precedes vertex $z_2$ on $\overrightarrow{Q}$ if $z_1$, $z_2$ occur on $\overrightarrow{Q}$ in this order, and indicate this relationship by $z_1\prec z_2$.

Let $P=x\overrightarrow{P}y$ be a path. A vine on $P$ is a set
$$
\{L_i=x_i\overrightarrow{L}_iy_i: 1\le i\le m\}
$$
of internally-disjoint paths such that\\

$(a)$ $V(L_i)\cap V(P)=\{x_i,y_i\} \ \ (i=1,...,m)$,

$(b)$ $x=x_1\prec x_2\prec y_1\preceq x_3\prec y_2\preceq x_4\prec ...\preceq x_m\prec y_{m-1}\prec y_m=y$ on $P$.\\

\noindent\textbf{The Vine Lemma} {[2]}. Let $G$ be a $k$-connected graph and $P$ a path in $G$. Then there are $k-1$ pairwise-disjoint vines on $P$.\\

The next three lemmas are crucial for the proof of Theorem 1.\\

\noindent\textbf{Lemma 1}. Let $G$ be a connected graph and $P=x\overrightarrow{P}y$ a longest path in $G$.

$(i)$ If $xz,yz^-\in E(G)$ for some $z\in V(x^+\overrightarrow{P}y)$ then $c=p$, that is $G$ is hamiltonian.

$(ii)$ If $d(x)+d(y)\ge p$ then $c=p$.

$(iii)$ Let $yz_1, xz_2\in E(G)$ for some $z_1,z_2 \in V(P)$ with $x\prec z_1\prec z_2 \prec y$ and $|z_1\overrightarrow{P}z_2|\ge 3$. If $xz,yz\not\in E(G)$ for each $z\in V(z_1^+\overrightarrow{P}z_2^-)$ and $d(x)+d(y)\ge p+3-|z_1\overrightarrow{P}z_2|$ then $c=p$. \\

\noindent\textbf{Lemma 2}. Let $G$ be a graph and $\{L_1,L_2,...,L_m\}$ be a vine on a longest cycle of $G$. Then
$$
c\ge\frac{2p-10}{m+1}+4.
$$

\noindent\textbf{Lemma 3}. Let $G$ be a connected graph and $\{L_1,L_2,...,L_m\}$ be a vine on a longest path $P=x\overrightarrow{P}y$ of $G$. Then either $c=p$ or $c\ge d(x)+d(y)+m-2$.\\

\section{Proofs}

\noindent\textbf{Proof of Lemma 1}. $(i)$ Let $xz,yz^-\in E(G)$ for some $z\in V(x^+\overrightarrow{P}y)$. Then $c\ge|xz\overrightarrow{P}yz^-\overleftarrow{P}x|=p$. If $V(G)=V(P)$ then clearly $c=p$. Otherwise, recalling that $G$ is connected, we can form a path longer that $P$, a contradiction.

$(ii)$ Let $d(x)+d(y)\ge p$. If $xz,yz^-\in E(G)$ for some $z\in V(x^+\overrightarrow{P}y)$ then we can argue as in $(i)$. Otherwise $N(x)\cap N^+(y)=\emptyset$. Observing also that $x\not\in N(x)\cup N^+(y)$, we get
$$
p\ge |N(x)|+|N^+(y)|+|\{x_1\}|
$$
$$
=|N(x)|+|N(y)|+1=d(x)+d(y)+1,
$$
contradicting the hypothesis.

$(iii)$ Assume the contrary, that is $c\le p-1$. Then by $(i)$, $N(x)\cap N^+(y)=\emptyset$. Clearly, $x\not\in N(x)\cup N^+(y)$. Further, by the hypothesis,
$$
V(z_1^{+2}\overrightarrow{P}z_2^-)\cap (N(x)\cup N^+(y))=\emptyset,
$$
implying that
$$
p\ge |\{x\}|+|N(x)|+|N^+(y)|+|V(z_1^{+2}\overrightarrow{P}z_2^-)|
$$
$$
=d(x)+d(y)+|z_1\overrightarrow{P}z_2|-2,
$$
contradicting the hypothesis. Thus, $c=p$. Lemma 1 is proved.  \quad\quad   \qquad \rule{7pt}{6pt}\\

\noindent\textbf{Proof of Lemma 2}. Let $P=x\overrightarrow{P}y$ be a longest path in $G$. Put

$$
L_i=x_i\overrightarrow{L}_iy_i \ \ (i=1,...,m), \ \ A_1=x_1\overrightarrow{P}x_2, \ \ A_m=y_{m-1}\overrightarrow{P}y_m,
$$
$$
A_i=y_{i-1}\overrightarrow{P}x_{i+1} \ \ (i=2,3,...,m-1), \ \  B_i=x_{i+1}\overrightarrow{P}y_i  \  \  (i=1,...,m-1),
$$
$$
|A_i|-1=a_i \  \ (i=1,...,m),\  \ |B_i|-1=b_i \  \ (i=1,...,m-1).
$$

By combining appropriate $L_i, A_i, B_i$, we can form the following cycles:
$$
Q_1=\bigcup_{i=1}^mA_i\cup\bigcup_{i=1}^mL_i,
$$
$$
Q_2=\bigcup_{i=1}^{m-1}A_i\cup B_{m-1}\cup\bigcup_{i=1}^{m-1}L_i,
$$
$$
Q_3=\bigcup_{i=2}^{m}A_i\cup B_{1}\cup\bigcup_{i=2}^{m}L_i,
$$
$$
R_i=B_i\cup A_{i+1}\cup B_{i+1}\cup L_{i+1} \ \  (i=1,...,m-2).
$$
Since $|L_i|\ge2$ $(i=1,...,m)$, we have
$$
c\ge|Q_1|=\sum_{i=1}^ma_i+\sum_{i=1}^m(|L_i|-1)\ge\sum_{i=1}^ma_i+m,
$$
$$
c\ge|Q_2|=b_{m-1}+\sum_{i=1}^{m-1}a_i+\sum_{i=1}^{m-1}(|L_i|-1)\ge b_{m-1}+\sum_{i=1}^{m-1}a_i+m-1,
$$
$$
c\ge|Q_3|=b_{1}+\sum_{i=2}^{m}a_i+\sum_{i=2}^{m}(|L_i|-1)\ge b_{1}+\sum_{i=2}^{m}a_i+m-1,
$$
$$
c\ge|R_i|=b_i+a_{i+1}+b_{i+1}+|L_{i+1}|-1
$$
$$
\ge b_i+a_{i+1}+b_{i+1}+1 \ \  (i=1,...,m-2).
$$
By summing, we get
$$
(m+1)c\ge\left(2\sum_{i=1}^{m}a_i+2\sum_{i=1}^{m-1}b_i\right)+2\sum_{i=2}^{m-1}a_i+4m-4
$$
$$
\ge2\left(\sum_{i=1}^ma_i+\sum_{i=1}^{m-1}b_i+1\right)+4m-6=2p+4m-6,
$$
implying that
$$
c\ge\frac{2p-10}{m+1}+4.
$$
Lemma 2 is proved.  \quad\quad\quad\quad\quad   \qquad \rule{7pt}{6pt}\\

\noindent\textbf{Proof of Lemma 3}. If $m=1$ then $xy\in E(G)$ and by Lemma 1$(i)$, $c=p$. Let $m\ge2$. Put $L_i=x_i\overrightarrow{L}_iy_i$ \ $(i=1,...,m)$ and let
$$
A_i, \ B_i, \ a_i, \ b_i, \ Q_i
$$
be as defined in the proof of Lemma 2.\\

\textbf{Case 1}. $m=2$.

Assume without loss of generality that $L_1$ and $L_2$ are chosen so as to minimize $b_1$. This means that $N(x)\cup N(y)\subseteq V(A_1\cup A_2)$. By Lemma 1$(iii)$, either $c=p$ or $d(x)+d(y)\le p+2-|z_1\overrightarrow{P}z_2|=p+1-b_1$. If $c=p$ then we are done. Let $d(x)+d(y)\le p+1-b_1$, that is $p\ge d(x)+d(y)+b_1-1$. Then $p=a_1+a_2+b_1+1\ge d(x)+d(y)+b_1-1$, implying that
$$
c\ge|Q_1|=a_1+a_2+2\ge d(x)+d(y)=d(x)+d(y)+m-2.
$$

\textbf{Case 2}. $m=3$.

Let $xz_1,yz_2\in E(G)$ for some $z_1,z_2\in V(P)$. If $z_2\prec z_1$ then $\{xz_1,yz_2\}$ is a vine consisting of two paths (edges) and we can argue as in Case 1. Otherwise we have
$$
N(x)\subseteq V(A_1\cup A_2), \ \  N(y)\subseteq V(A_2\cup A_3)
$$
and $z_1\preceq z_2$ for each $z_1\in N(x)$ and $z_2\in N(y)$. Therefore, $a_1+a_2+a_3\ge d(x_1)+d(x_2)-2$ and
$$
c\ge|Q_1|=a_1+a_2+a_3+3
$$
$$
\ge d(x_1)+d(x_2)+1=d(x_1)+d(x_2)+m-2.
$$

\textbf{Case 3.} $m\ge4$.

Choose $\{L_1,...,L_m\}$  so as to minimize $m$. Then clearly
$$
N(x)\subseteq V(A_1\cup A_2), \ \  N(y)\subseteq V(A_{m-1}\cup A_{m})
$$
and $z_1\prec z_2$ for each $z_1\in N(x)$ and $z_2\in N(y)$. Observing also that
$$
a_1+a_2\ge d(x)-1, \ \ a_{m-1}+a_m\ge d(y)-1,
$$
we get
$$
c\ge|Q_1|=\sum_{i=1}^ma_i+m=(a_1+a_2+a_{m-1}+a_m)+\sum_{i=3}^{m-3}a_i+m
$$
$$
=d(x)+d(y)-2+\sum_{i=3}^{m-3}a_i+m\ge d(x)+d(y)=m-2.
$$
Lemma 3 is proved.      \quad \quad \quad  \quad \rule{7pt}{6pt}\\

\noindent\textbf{Proof of Theorem 1}. Let $P=x\overrightarrow{P}y$ be  a longest path in $G$. \\

\textbf{Case 1}. $p\le 2\delta$.

If $xy\in E(G)$ then by Lemma 1$(i)$, $c=p$. Let $xy\not\in E(G)$. Then $d(x)+d(y)\ge2\delta\ge p$ and by Lemma 1$(ii)$, $c=p$. \\

\textbf{Case 2}. $2\delta+1\le p\le 3\delta-2$.

If $c\ge 3\delta-3$ then $c\ge p+2-3=p-1$. Next, if $\kappa=2$ and $p\ge 3\delta-1$ then $p\ge 3\delta-1\ge p+1$, a contradiction. By Theorem C, $c\ge p-1$. \\

\textbf{Case 3}. $p\ge 3\delta-1$.

Since $G$ is 2-connected, there is a vine $\{L_1,...,L_m\}$ on $P$. By Lemma 3, $m\le c-d(x)-d(y)+2\le c-2\delta+2$. Using Lemma 2, we get
$$
c\ge \frac{2p-10}{m+1}+4\ge\frac{2p-10}{c-2\delta+3}+4,
$$
implying that
$$
c\ge \sqrt{2p-10+\left(\delta-\frac{7}{2}\right)}+\delta+\frac{1}{2}.
$$
Theorem 1 is proved.     \quad\quad\quad\quad   \quad \rule{7pt}{6pt}

\noindent Institute for Informatics and Automation Problems\\ National Academy of Sciences\\
P. Sevak 1, Yerevan 0014, Armenia\\ E-mail: zhora@ipia.sci.am

\end{document}